\documentclass{amsart}

\newcommand{\Z}{\ensuremath{\mathbf Z}}

\newcommand{\N}{\ensuremath{ \mathbf N }}

\newcommand{\R}{\ensuremath{ \mathbf{R} }}

\newcommand{\mca}{\ensuremath{ \mathcal{A} }}
\newcommand{\mcb}{\ensuremath{ \mathcal{B} }}

\newtheorem{theorem}{Theorem}
\newtheorem{lemma}{Lemma}
\newtheorem{corollary}{Corollary}
\newcommand{\bt}{\begin{theorem}}
\newcommand{\et}{\end{theorem}}
\newcommand{\bl}{\begin{lemma}}
\newcommand{\el}{\end{lemma}}
\newcommand{\bc}{\begin{corollary}}
\newcommand{\ec}{\end{corollary}}

\newtheorem{problem}{Problem} 
\newcommand{\bprob}{\begin{problem}}
\newcommand{\eprob}{\end{problem}}

\newcommand{\bpf}{\begin{proof}}
\newcommand{\epf}{\end{proof}}

\newcommand{\beq}{\begin{equation}}
\newcommand{\eeq}{\end{equation}}
\newcommand{\benum}{\begin{enumerate}}
\newcommand{\eenum}{\end{enumerate}}
\newcommand{\card}{\ensuremath{\text{card}}}
\newcommand{\diam}{\ensuremath{\text{diam}}}

\title{Linear forms and complementing sets of integers}
\author{Melvyn B. Nathanson}



\thanks{This work was supported
in part by grants from the NSA Mathematical Sciences Program and
the PSC-CUNY Research Award Program}
\address{Department of Mathematics, Lehman College (CUNY), 
Bronx, NY 10468}
\curraddr{School of Mathematics, Institute for Advanced Study,
Princeton, NJ 08540}
\email{melvyn.nathanson@lehman.cuny.edu}

\date{\today}

\begin{document}

\maketitle

\begin{abstract}
Let $\varphi(x_1,\ldots,x_h,y) = u_1x_1 + \cdots  + u_hx_h+vy$ be a  linear form with nonzero integer coefficients $u_1,\ldots,  u_h, v.$  
Let $\mathcal{A} = (A_1,\ldots, A_h)$ be an $h$-tuple of finite sets of integers and let $B$ be an infinite set of integers.  Define the representation function associated to the form $\varphi$ and the sets  \mca\ and $B$  as follows:
\[
R^{(\varphi)}_{\mca,B}(n) = 
 \text{card}\left( \left\{ (a_1,\ldots, a_h,b) \in A_1 \times \cdots \times A_h \times B: \varphi(a_1, \ldots , a_h,b ) = n \right\} \right).
\]
If this representation function is constant, then the set $B$ is periodic and the period of $B$ will be bounded in terms of the diameter of the finite set $\{ \varphi(a_1,\ldots,a_h,0): (a_1,\ldots, a_h) \in A_1 \times \cdots \times A_h\}.$
\end{abstract}

\section{Representation functions for linear forms}
Let $h \geq 1$ and let 
\[
\psi(x_1,\ldots,x_h) = u_1x_1 + \cdots  + u_hx_h
\]
be a  linear form with nonzero integer coefficients $u_1,\ldots,  u_h.$  
Let 
\[
\mathcal{A} = (A_1,\ldots, A_h)
\]
be an $h$-tuple of sets of integers. 
The \emph{image} of $\psi$ with respect to $\mathcal{A}$ is the set
\[
\psi(\mathcal{A}) 
= \left\{ \psi(a_1, \ldots , a_h ) :  (a_1,\ldots, a_h) \in A_1 \times \cdots \times A_h  \right\}.
\]
Then $\psi(\mathcal{A}) \neq \emptyset$  if and only if $A_i \neq \emptyset$ for all $i=1,\ldots, h.$   For $\psi(\mathcal{A}) \neq \emptyset$, we define the \emph{diameter of $\mathcal{A}$ with respect to $\psi$} by
\[
D^{(\psi)}_{\mca} = \diam(\psi(\mathcal{A})) = \sup(\psi(\mathcal{A})) - \inf(\psi(\mathcal{A})).
\]
We have $D^{(\psi)}_{\mca}  > 0$ if and only if $|A_i| > 1$ for some $i.$

For every integer $n$, we define the \emph{representation function} associated to $\psi$ by 
\[
R^{(\psi)}_{\mathcal{A}}(n) = \text{card}\left( \left\{ (a_1,\ldots, a_h) \in A_1 \times \cdots \times A_h : \psi(a_1, \ldots , a_h ) = n \right\} \right).
\]
Then $n\in \psi(\mca)$ if and only if $R^{(\psi)}_{\mathcal{A}}(n) > 0.$

Let $\ell \geq 1$ and let 
\[
\omega(y_1,\ldots,y_{\ell}) = v_1y_1 + \cdots  +v_{\ell}y_{\ell}
\]
be another linear form with nonzero integer coefficients $v_1,\ldots,  v_{\ell}.$  Consider the linear form
\[
\varphi(x_1,\ldots, x_h, y_1,\ldots,y_{\ell}) = \psi(x_1,\ldots,x_h) + \omega(y_1,\ldots,y_{\ell}).
\]
Let $\mathcal{A} = (A_1,\ldots, A_h)$ be an $h$-tuple of sets of integers and let $\mathcal{B} = (B_1,\ldots, B_{\ell})$ be an $\ell$-tuple of sets of integers.   The \emph{image} of $\varphi$ with respect to $(\mca,\mcb)$ is the set
\begin{align*}
\varphi(\mca,\mcb)  = & \psi(\mathcal{A}) + \omega(\mcb) \\
= & \left\{ \psi(a_1, \ldots , a_h ) + \omega(b_1,\ldots, b_{\ell}) : 
(a_1,\ldots, a_h) \in A_1 \times \cdots \times A_h   \right. \\
& \left. \text{ and } (b_1,\ldots, b_{\ell}) \in B_1 \times \cdots \times B_{\ell} \right\}.
\end{align*}
We define the \emph{representation function associated to $\varphi$, \mca, and \mcb}  by 
\begin{align*}
R^{(\varphi)}_{\mca,\mcb}(n) = 
& \text{card}\left( \left\{ (a_1,\ldots, a_h,b_1,\ldots, b_{\ell}) \in A_1 \times \cdots \times A_h \times B_1 \times \cdots \times B_{\ell}: \right. \right. \\
& \left. \left. \varphi(a_1, \ldots , a_h,b_1,\ldots, b_{\ell} ) = n \right\} \right).
\end{align*}

For every positive integer $m$, we define the \emph{modular representation function} associated to $\varphi$ by 
\begin{align*}
R^{(\varphi)}_{\mca,\mcb;m}(n) = 
& \text{card}\left( \left\{ (a_1,\ldots, a_h,b_1,\ldots, b_{\ell}) \in A_1 \times \cdots \times A_h \times B_1 \times \cdots \times B_{\ell}: \right. \right. \\
& \left. \left. \varphi(a_1, \ldots , a_h,b_1,\ldots, b_{\ell} )  \equiv n \pmod{m}  \right\} \right).
\end{align*}

If $\ell =1$ and $\mcb = (B),$ then we write $\varphi(\mca,\mcb) = \varphi(\mca,B),$ $R^{(\varphi)}_{\mca,\mcb}(n) = R^{(\varphi)}_{\mca,B}(n),$
and $R^{(\varphi)}_{\mca,\mcb;m}(n) = R^{(\varphi)}_{\mca,B;m}(n).$

\emph{Notation.}  
Let \Z\ and $\N_0$ denote the set of integers and the set of nonnegative integers, respectively.  
We denote the cardinality of the set $S$ by $|S|$  or by $\card(S).$
We denote  the integer part of the real number $x$  by $[x]$.

\section{complementing sets}

A classical problem in additive number theory is the study of complementing pairs of sets of integers, that is, pairs $(A,B)$ such that every integer has a unique representation in the form $a+b,$ with $a\in A$ and $b \in B.$   This is the case $h=1,$ $\psi(x) = x,$  $\omega(y) = y,$ and $\varphi(x,y) = x+y$ of the general problem of representations of integers by linear forms.  There are many beautiful results and open problems about complementing sets for the integers.  For example, if $A$ is a finite set of integers and if $B$ is an infinite set of integers such that the pair $(A,B)$ is complementing, then $B$ must be a periodic set, that is, a union of congruence classes modulo $m$ for some positive integer $m$ (Newman~\cite{newm77}).  There are also upper and lower bounds on the period $m$ as a function of the diameter of the set $A$ (Biro~\cite{biro05}, Kolountzakis~\cite{kolo03}, Ruzsa~\cite[Appendix]{tijd06}, Steinberger~\cite{stei07}).  
In general,  it is known that every pair $(A,B)$ of complementing sets with $A$ finite must satisfy a certain cyclotomy condition, but it is an open problem to determine if a finite set $A$ of integers has a complement.  

Complementing pairs have also been studied for sets of lattice points (Hansen~\cite{hans69}, Nathanson~\cite{nath72c}, Niven~\cite{nive71}).  
If $(A,B)$ is a pair of sets of lattice points such that $A$ is finite and every lattice point has a unique representation in the form $a+b$ with $a\in A$ and $b\in B,$ then it is an open problem to determine if the set $B$ must be periodic.

The object of this paper is to begin the study of complementing sets of integers with respect to an arbitrary linear form $\varphi(x_1,\ldots,x_h,y_1,\ldots, y_{\ell}).$ 
Let \mca\ be an $h$-tuple of sets of integers and \mcb\ an $\ell$-tuple of sets of integers.  The pair $(\mca,\mcb)$ is called \emph{complementing} with respect to $\varphi$ 
if $R^{(\varphi)}_{\mca,\mcb}(n)=1$ for all $n \in \mathbf{Z},$   that is, if every integer $n$ has a unique representation in the form $n=\psi(a_1,\ldots, a_h) + \omega(b_1,\ldots, b_{\ell}) ,$ where $a_i \in A_i$ for $i=1,\ldots, h$ and $b_j \in B_j$ for $j=1,\ldots, \ell.$   The pair $(\mca,\mcb)$  is called \emph{$t$-complementing} with respect to $\varphi$  if $R^{(\varphi)}_{\mca,\mcb}(n)=t$ for all $n\in \mathbf{Z}.$  
The pair $(\mca,\mcb)$   is called \emph{$t$-complementing  modulo $m$} with respect to $\varphi$  if $R^{(\varphi)}_{\mca,\mcb;m}(\ell)=t$ for all $\ell \in \{0,1,\ldots, m-1\}.$

The pair $(\mca,\mcb)$  is called \emph{periodic with respect to $\varphi$}  if the representation function $R^{(\varphi)}_{\mca,\mcb}$ is  periodic, that is, if there is a positive integer $m$ such that 
$R^{(\varphi)}_{\mca,\mcb}(n+m) =  R^{(\varphi)}_{\mca,\mcb}(n)$ for all integers $n.$  
The pair $(\mca,\mcb)$  is called \emph{eventually periodic with respect to $\varphi$} if the representation function $R^{(\varphi)}_{\mca,\mcb}$ is eventually periodic, that is,  if there exist integers $m \geq 1$ and $n_0$ such that 
$R^{(\varphi)}_{\mca,\mcb}(n+m) =  R^{(\varphi)}_{\mca,\mcb}(n)$ for all integers $n\geq n_0.$  

We consider the case $\ell = 1.$   Suppose that $\varphi(x_1,\ldots,x_h,y) = \psi(x_1,\ldots,x_h) + vy$ is a linear form with nonzero integer coefficients, and that $\mca$ is an $h$-tuple of finite sets of integers and $B$ is a set of integers such that the pair $(\mca,B)$  is $t$-complementing with respect to $\varphi.$ 
We shall prove that the set $B$ is periodic, and obtain an upper bound for the period of $B$ in terms of the diameter $D^{\psi}_{\mca}$ of the finite set $\psi(\mca).$  We also obtain a cyclotomic condition related to $t$-complementing sets modulo $m,$ and describe a compactness argument that allows us to solve an inverse problem related to representation functions associated with linear forms.

\section{Linear forms and periodicity}

\bt
Let $h \geq 1$ and let 
\[
\varphi(x_1,\ldots,x_h,y) = u_1x_1 + \cdots +  u_hx_h + vy
\]
be a  linear form with nonzero integer coefficients $u_1,\ldots,  u_h, v.$ 
Let  $\mathcal{A} = (A_1,\ldots, ,A_h)$ be an $h$-tuple of nonempty finite sets of integers, and let $B$ be an infinite set of integers.  If $(\mca,B)$ is $t$-complementing respect to $\varphi,$ then $B$ is periodic, that is, there is a positive integer $m$ such that $B$ is a union of congruence classes modulo $m$.
\et

\bpf
If $v < 0,$ then we replace $\varphi$ with $-\varphi.$  Thus, we can assume without loss of generality that $v \geq 1$.

If $|A_i|=1$ for all $i=1,\ldots, h,$ then the linear form $\varphi$ represents all integers if and only if $v=1$ and $B=\Z,$ and the Theorem holds with $m=1.$     Thus, we can also assume that $|A_i| > 1$ for at least one $i$.

Consider the linear form
\[
\psi(x_1,\ldots,x_h) = u_1x_1 + \cdots + u_hx_h.
\]
We have
\[
\varphi(a_1,\ldots, a_h, b) = \psi(a_1,\ldots, a_b) + vb
\]
for all $(a_1,\ldots, a_h) \in A_1\times \cdots \times A_h$ and $b\in B.$
Let $g_{\text{min}} = \min\left( \psi(A_1,\ldots, A_{h}) \right)$ and 
$g_{\text{max}} = \max\left( \psi(A_1,\ldots, A_{h}) \right)$.  
Since $|A_i| > 1$ for some $i \in \{1,2,\ldots, h\}$, it follows that $g_{\text{min}} < g_{\text{max}}$ and 
\[
D^{(\psi)}_{\mca} = \diam(\psi(A_1,\ldots, A_{h})) = g_{\text{max}} - g_{\text{min}}\geq 1.
\]
Let
\[
G_{\text{min}} = \left\{ ( a_{1}, \ldots, a_{h}) \in A_1\times\cdots\times  A_{h} : \psi( a_{1}, \ldots, a_{h}) = g_{\text{min}} \right\}
\]
and
\[
G_{\text{max}} = \left\{ ( a_{1}, \ldots, a_{h}) \in A_1\times\cdots\times  A_{h} : \psi( a_{1}, \ldots, a_{h}) = g_{\text{max}} \right\}.
\]
Then
\[
|G_{\text{min}}| = R^{(\psi)}_{\mca}(g_{\text{min}})  \geq 1
\]
and
\[
|G_{\text{max}}| = R^{(\psi)}_{\mca}(g_{\text{max}}) \geq 1.
\]
Let $\chi_{B}: \R \rightarrow \{0,1\}$ denote the characteristic function of the set $B,$ that is,
\[
\chi_{B}(x) = 
\begin{cases}
1 & \text{if $x \in B$} \\
0 & \text{if $x \notin B$.}
\end{cases}
\]
We have
\[
\varphi( a_{1},\ldots, a_h,b) = \psi( a_{1},\ldots, a_h) + vb= n
\]
 if and only if 
\[
b = \frac{n-\psi( a_{1}, \ldots, a_{h}) }{v} \in B.
\]
It follows that  
\[
R^{(\varphi)}_{\mca}(n)
 =  \sum_{( a_{1}, \ldots, a_{h}) \in A_1\times\cdots\times  A_{h}} \chi_{B}\left( \frac{n- \psi( a_{1}, \ldots, a_{h}) }{v}\right)
\]
for all $n \in \Z.$  
We can also write 
\begin{align*}
R^{(\varphi)}_{\mca}(n)
 =  &  \sum_{
\substack{ ( a_{1}, \ldots, a_{h}) \in A_1\times\cdots\times  A_{h}\\
( a_{1}, \ldots, a_{h}) \notin G_{\text{min}}}
} \chi_{B}\left( \frac{n- \psi( a_{1}, \ldots, a_{h}) }{v}\right) \\
& + |G_{\text{min}}| \chi_{B}\left( \frac{n- g_{\text{min}} }{v}\right).
\end{align*}
Replacing $n$ by $vn + g_{\text{min}}$, we obtain the identity
\begin{align*}
R^{(\varphi)}_{\mca}(v n + g_{\text{min}})
 =  &  \sum_{
\substack{ ( a_{1}, \ldots, a_{h}) \in A_1\times\cdots\times  A_{h}\\
( a_{1}, \ldots, a_{h}) \notin G_{\text{min}}}
} \chi_{B}\left( \frac{v n + g_{\text{min}}- \psi( a_{1}, \ldots, a_{h}) }{v}\right) \\
& + |G_{\text{min}}| \chi_{B}\left( n\right).
\end{align*}
Equivalently,
\begin{align*}
|G_{\text{min}}| \chi_{B}(n) 
= & R^{(\varphi)}_{\mca}(v n + g_{\text{min}}) \\
& -  \sum_{\substack{ ( a_{1}, \ldots, a_{h}) \in A_1\times\cdots\times  A_{h}\\
( a_{1}, \ldots, a_{h}) \notin G_{\text{min}}}}
\chi_{B}\left( n - \frac{\psi( a_{1}, \ldots, a_{h}) - g_{\text{min}}}{v}\right).
\end{align*}
Since  $g_{\text{min}} < \psi( a_{1}, \ldots, a_{h}) \leq g_{\text{max}}$ for all $h$-tuples $( a_{1}, \ldots, a_{h}) \notin G_{\text{min}},$ it follows that 
\[
0 < \frac{1}{v} \leq \frac{\psi( a_{1}, \ldots, a_{h}) - g_{\text{min}}}{u_h} \leq 
\frac{g_{\text{max}} - g_{\text{min}}}{u_h}.
\]

Similarly, replacing $n$ by $v n + g_{\text{max}}$, we obtain the identity
\begin{align*}
|G_{\text{max}}| \chi_{B}(n) 
= & R^{(\varphi)}_{\mca}(v n + g_{\text{max}}) \\
& -  \sum_{
\substack{ ( a_{1}, \ldots, a_{h}) \in A_1\times\cdots\times  A_{h}\\
( a_{1}, \ldots, a_{h}) \notin G_{\text{min}}}} 
\chi_{B}\left( n + \frac{g_{\text{max}}  - \psi( a_{1}, \ldots, a_{h}) }{v}\right).
\end{align*}
Since  $g_{\text{min}} \leq \psi( a_{1}, \ldots, a_{h}) < g_{\text{max}}$ for $( a_{1}, \ldots, a_{h}) \notin G_{\text{max}},$ it follows that 
\[
0 < \frac{1}{v} \leq \frac{ g_{\text{max}} - \psi( a_{1}, \ldots, a_{h}) }{v} \leq 
\frac{g_{\text{max}} - g_{\text{min}}}{v}.
\]

We define the nonnegative integer
\beq  \label{CS:diam}
d =  \left[\frac{\diam(\psi(A_1,\ldots,A_{h}))}{v}\right] 
=  \left[\frac{g_{\text{max}} - g_{\text{min}}}{v}\right].
\eeq
If the pair $(\mca, B)$ is $t$-complementing with respect to $\varphi,$ then $R^{(\varphi)}_{\mca,B}(n)=t$ for all $n \in \Z,$ and so 
\[
|G_{\text{min}}| \chi_{B}(n) 
= t  -   \sum_{
\substack{ ( a_{1}, \ldots, a_{h}) \in A_1\times\cdots\times  A_{h}\\
( a_{1}, \ldots, a_{h}) \notin G_{\text{min}}}} 
 \chi_{B}\left( n - \frac{\psi( a_{1}, \ldots, a_{h}) - g_{\text{min}}}{v}\right)
\]
and
\[
|G_{\text{max}}| \chi_{B}(n) 
= t -   \sum_{
\substack{ ( a_{1}, \ldots, a_{h}) \in A_1\times\cdots\times  A_{h}\\
( a_{1}, \ldots, a_{h}) \notin G_{\text{min}}}} 
\chi_{B}\left( n + \frac{g_{\text{max}}  - \psi( a_{1}, \ldots, a_{h}) }{v}\right)
\]
These formulae allow us to compute the characteristic function $\chi_{B}$ recursively for all integers if we know the value of $\chi_{B}$ for any $d$ consecutive integers.

Consider the $d$-tuple 
\[
\mathcal{B}(j) = (\chi_{B}(j), \chi_{B}(j+1),\ldots, \chi_{B}(j+d-1)) \in \{ 0,1 \}^d.
\]
Since there only $2^d$ binary sequences of length $d$, it follows from the pigeonhole principle that there are integers $j_1, j_2$ such that $0 \leq j_1 < j_2 \leq 2^d$ and $\mathcal{B}(j_1) = \mathcal{B}(j_2).$  Let $m = j_2 - j_1.$   Then 
\[
1 \leq m \leq 2^d
\]
and $\chi_{B}(n) = \chi_{B}(n+m)$ for $n=j_1,\ldots, j_1+d-1.$
The recursion formulae imply that $\chi_{B}(n) = \chi_{B}(n+m)$ for all integers $n$.  
This completes the proof.
\epf

\section{Linear forms and cyclotomy}

\bt
Let $h \geq 1$ and let 
\[
\psi(x_1,\ldots,x_h,y) = u_1x_1 + \cdots  + u_hx_h
\]
be a  linear form with nonzero integer coefficients $u_1,\ldots, u_h,.$  
Let   $\mathcal{A} = (A_1,\ldots, A_h)$ be an $h$-tuple of nonempty finite sets of integers.   Consider the modular representation function 
\[
R^{(\psi)}_{\mathcal{A},m}(n) = \card\left( \{ (a_1,\ldots, a_h) \in A_1 \times \cdots \times A_h: \psi(a_1,\ldots, a_h) \equiv n \pmod{m} \}\right).
\]
and the generating functions
\[
F_{A_i}(z) = \sum_{a_i\in A_i} z^{a_i} \qquad\text{for $i=1,\ldots, h$}.
\]
For $m \geq 1$, define the polynomial 
\[
\Lambda_m(z) = 1 + z + z^2 + \cdots + z^{m-1}.
\]
The $h$-tuple $\mathcal{A}$ is $t$-complementing modulo $m$ with respect to $\psi$  if and only if there exists a nonnegative integer $L$ such that 
\beq  \label{CS:condition}
z^L F_{A_1}(z^{u_1}) \cdots F_{A_h}(z^{u_h}) \equiv t\Lambda_m(z) \pmod{z^m-1}.
\eeq
\et

\bpf
The generating functions $F_{A_i}(z)$ are nonzero Laurent polynomials for $i=1,\ldots, h.$   The function 
\[
F(z) = F_{A_1}(z^{u_1}) \cdots F_{A_h}(z^{u_h})
\]
is also a nonzero Laurent polynomial with integer coefficients.    Choose a nonnegative integer $L$ such that  $z^{L}F(z)$ is a polynomial. 

The sets $A_1,\ldots, A_h$ are finite, and so $\psi(\mathcal{A})$ is finite.  We have $ R^{(\psi)}_{\mathcal{A}}(n) \geq 1$ if and and only if $n \in \psi(\mathcal{A})$.  For $\ell = 0,1,\ldots, m-1,$ we consider the finite set  
\[
\mathcal{I}_{\ell} = \{ i \in \Z :  R^{(\psi)}_{\mathcal{A}}(\ell + im) \geq 1\}.
\]
Since $F_{A_i}(z^{u_i}) = \sum_{a_i\in A_i}z^{u_ia_i}$ for $i=1,\ldots, h,$ it follows that 
\begin{align*}
F(z) & = F_{A_1}(z^{u_1}) \cdots F_{A_h}(z^{u_h}) \\
& =\sum_{a_1 \in A_1 } \cdots \sum_{a_h \in A_h } z^{u_1a_1+\cdots + u_ha_h} \\
& =\sum_{a_1 \in A_1 } \cdots \sum_{a_h \in A_h } z^{\psi(a_1,\ldots, a_h)} \\
& = \sum_{n \in \psi(\mathcal{A})} R^{(\psi)}_{\mathcal{A}}(n)z^n \\
& = \sum_{\ell =0}^{m-1}  \sum_{\substack{n\in \psi(\mathcal{A}) \\ n\equiv \ell\pmod{m}} }R^{(\psi)}_{\mathcal{A}}(n) z^n \\
& = \sum_{\ell =0}^{m-1}  \sum_{i\in \mathcal{I}_{\ell}}  R^{(\psi)}_{\mathcal{A}}(\ell+im)z^{\ell + im}.
\end{align*}
Since 
\[
z^{L}F(z) =  \sum_{\ell =0}^{m-1}  \sum_{i\in \mathcal{I}_{\ell}}  R^{(\psi)}_{\mathcal{A}}(\ell+im)z^{\ell + L + im}
\]
is a polynomial, it follows that $\ell + L+im \geq 0$ for all $\ell \in \{0,1,\ldots, m-1\}$ and $i \in \mathcal{I}_{\ell}$.  Applying the division algorithm for integers, we can write
\[
\ell + L = \alpha(\ell)+ \beta(\ell) m
\]
where $0 \leq \alpha(\ell) \leq m-1$ for $\ell = 0,1,\ldots, m-1.$  
Moreover, if $\ell \not\equiv \ell' \pmod{m},$ then $\alpha(\ell) \neq \alpha(\ell')$ and so
\[
\{\alpha(0), \alpha(1), \ldots, \alpha(m-1)\} = \{ 0,1,\ldots, m-1\}.
\]
Equivalently,
\[
\sum_{\ell=0}^{m-1}z^{\alpha(\ell)} = \sum_{\ell=0}^{m-1}z^{\ell} = \Lambda_m(z).
\]
If $i \in \mathcal{I}_{\ell},$ then the inequality 
\[
\ell + L+im = \alpha(\ell)+ (\beta(\ell) + i) m\geq 0
\]
implies that $\beta(\ell) + i \geq 0$.  Therefore, for each $\ell \in \{0,1,\ldots, m-1\}$ there is a polynomial $p_{\ell}(z)$ with integral coefficients such that 
\begin{align*}
\sum_{i\in \mathcal{I}_{\ell}}  R^{(\psi)}_{\mathcal{A}}(\ell+im)z^{\ell + L + im}  
&   =   \sum_{i\in \mathcal{I}_{\ell}}  R^{(\psi)}_{\mathcal{A}}(\ell+im)z^{\alpha(\ell)+ (\beta(\ell) + i) m}  \\
&   =   \sum_{i\in \mathcal{I}_{\ell}}  R^{(\psi)}_{\mathcal{A}}(\ell+im)z^{\alpha(\ell)}\left( 1 + ( z^m -1)  \right)^{\beta(\ell) + i}\\
&   =   \sum_{i\in \mathcal{I}_{\ell}}  R^{(\psi)}_{\mathcal{A}}(\ell+im)z^{\alpha(\ell)}   +  ( z^m -1) p_{\ell}(z) \\
&   =   R^{(\psi)}_{\mathcal{A},m}(\ell)z^{\alpha(\ell)}   +  ( z^m -1) p_{\ell}(z).
\end{align*}
It follows that 
\begin{align*}
z^{L}F(z)
&  =  \sum_{\ell =0}^{m-1}  \sum_{i\in \Z}  R^{(\psi)}_{\mathcal{A}}(\ell+im)z^{\ell + L + im}  \\
&   =  \sum_{\ell =0}^{m-1}  R^{(\psi)}_{\mathcal{A},m}(\ell)z^{\alpha(\ell)}   +   ( z^m -1)\sum_{\ell =0}^{m-1} p_{\ell}(z) \\
& = r_L(z) + (z^m - 1) q_L(z), 
\end{align*}
where 
\[
q_L(z) = \sum_{\ell =0}^{m-1} p_{\ell}(z) 
\]
and
\[
r_L(z) = \sum_{\ell =0}^{m-1}  R^{(\psi)}_{\mathcal{A},m}(\ell)z^{\alpha(\ell)}
\]
is a polynomial of degree at most $m-1$.
By the division algorithm for polynomials, this representation of $z^LF(z)$  is unique.  

Suppose that $\mathcal{A} = (A_1,\ldots, A_h)$ is a $t$-complementing $h$-tuple modulo $m$.  Then $R_{\mathcal{A},m}(\ell)=t$ for all $\ell$, and 
\[
r_L(z) = \sum_{\ell =0}^{m-1} t z^{\alpha(\ell)} = t\Lambda_m(z).
\]
Therefore, 
\[
z^{L}F(z) =  t\Lambda_m(z) + (z^m-1)q_L(z)
\]
and condition~\eqref{CS:condition} is satisfied.

Conversely, suppose that the generating functions $F_{A_1}(z),\ldots, F_{A_h}(z)$ satisfy condition~\eqref{CS:condition} for some nonnegative integer $L$.  By the uniqueness of the polynomial division algorithm, we have
\[
\sum_{\ell = 0}^{m-1} t z^{\ell} = t\Lambda_m(z) = r_L(z) = \sum_{\ell =0}^{m-1}  R^{(\psi)}_{\mathcal{A},m}(\ell) z^{\alpha(\ell)}.
\] 
Since 
\[
\{\alpha(0), \alpha(1), \ldots, \alpha(m-1)\} = \{ 0,1,\ldots, m-1\},
\]
it follows that $R^{(\psi)}_{\mathcal{A},m}(\ell) = t$ for all $\ell \in \{ 0,1,\ldots, m-1\},$ and so $\mathcal{A} = (A_1,\ldots, A_h)$ is a $t$-complementing $h$-tuple modulo $m$.  This completes the proof.
\epf

\section{An inverse problem for linear forms}
There are several inverse problems for representation functions associated to linear forms.  For example, let $\varphi(x_1,\ldots, x_h,y)$ be a form in $h+1$ variables and let $f:\Z \rightarrow \N_0 \cup \{\infty\}$ be a function.  If $\mca = (A_1,\ldots, A_{h})$ is an  $h$-tuple of sets of integers, does there exist a set $B$ such that the pair $(\mca, B)$  satisfies $R^{(\varphi)}_{\mca,B}(n) = f(n)$ for all $n \in \Z$?  In this section we use a compactness argument to obtain a result  in the case that $\mca = (A_1,\ldots, A_{h})$ is an  $h$-tuple of finite sets.

\bt  \label{CS:theorem:compact}
Let $h \geq 1$ and let 
\[
\varphi(x_1,\ldots,x_h,y) = u_1x_1 + \cdots + u_hx_h + u_hx_h+vy
\]
be a  linear form with nonzero integer coefficients $u_1,\ldots,u_h, v.$  
Let  
$\mca = (A_1,\ldots, A_{h})$
be an $h$-tuple of nonempty finite sets of integers.   Let $f:\Z \rightarrow \N_0 \cup \{\infty\}$ be a function.  
Suppose that there is a strictly increasing sequence $\{L_N\}_{N=1}^{\infty}$ of positive integers with the property that, for every $N \geq 1$, there exists a set $B_N$ of integers that satisfies 
\[
R^{(\varphi)}_{\mca,B_N}(n) = f(n) \qquad\text{for $|n| \leq L_N.$}
\]
Then there exists a set $B$ such that 
\[
R^{(\varphi)}_{\mca,B}(n) = f(n) \qquad\text{for all $n \in \Z.$}
\]
\et

\bpf
Since $L_N \geq N$ for all $N \geq 1,$ we can assume without loss of generality that $L_N = N.$  Consider the linear form
\[
\psi(x_1,\ldots,x_{h}) = u_1x_1 + \cdots + u_{h}x_{h} 
\]
Then
\[
\varphi(a_1,\ldots, a_h,b) = \psi(a_1,\ldots,a_{}) + vb
\]
for all integers $a_1,\ldots, a_h,b.$  
Moreover, since the sets $A_1,\ldots, A_{h}$ are finite, there is  a positive integer $g^{\ast}$ such that $\psi(\mca) \subseteq [-g^{\ast},g^{\ast}].$  
If $(a_1,\ldots,a_h) \in A_1 \times \cdots \times A_h$, if $b\in \Z$, and if $\varphi(a_1,\ldots, a_h,b)=  n \in [-N,N],$ then 
\[
v |b| = |n-\psi(a_1,\ldots, a_h) | \leq |n|+|\psi(a_1,\ldots, a_h) | \leq  N+g^{\ast}.
\]
Replacing the set $B_N$ with 
$B_N \cap [-(N+g^{\ast})/v, (N+g^{\ast})/v],$ 
we can assume without loss of generality that 
$B_N \subseteq  [-(N+g^{\ast})/v, (N+g^{\ast})/v]$ 
for all $N \geq 1.$

We shall construct inductively an increasing sequence of finite sets  
$B'_{1} \subseteq B'_{2}  \subseteq \cdots$ with the following properties:
\benum
\item     \label{CS:prop1}
For every positive integer $i$ and every integer $n \in [-i,i]$ we have $R^{(\varphi)}_{\mca,B'_i}(n) = f(n)$. 
\item        \label{CS:prop2}
For every positive integer $i$ there is a strictly increasing sequence $\left\{ N^{(i)}_j\right\}_{j=1}^{\infty}$ such that $i \leq N^{(i)}_1$ and $B'_i \subseteq B_{ N^{(i)}_j}$ for all $j \geq 1.$
\eenum

We begin by constructing the set $B'_{1}.$  
If $(a_1,\ldots,a_h) \in A_1 \times \cdots \times A_h$, if $b\in \Z$, and if $\varphi(a_1,\ldots, a_h,b) \in [-1,1],$ then 
\[
|b|\leq  \frac{1+g^{\ast}}{v}.
\]
For all $N \geq 1$ we have 
$R^{(\varphi)}_{\mca,B_N}(n) = f(n)$ for $|n| \leq N,$ 
and so $R^{(\varphi)}_{\mca,B_N}(n) = f(n)$ for $|n| \leq 1.$  
Let
\[
B^{(1)}_N = B_N \cap \left[-\frac{1+g^{\ast}}{v}, \frac{1+g^{\ast}}{v}\right]
\]
for $N \geq 1.$  Then $\left\{ B^{(1)}_N \right\}_{N=1}^{\infty}$ is an infinite sequence of subsets of the finite set 
$[-(1+g^{\ast})/v , (1+g^{\ast})/v ] \cap \Z$.  
By the pigeonhole principle, there is a strictly increasing sequence 
$\left\{N^{(1)}_j \right\}_{j=1}^{\infty}$ of positive integers and a set $B'_1$ such that $1 \leq N^{(1)}_1$ and 
\[
B'_1 = B^{(1)}_{N^{(1)}_j} \subseteq B_{N^{(1)}_j}
\]
for all $j\geq 1.$

Suppose that we have constructed an increasing sequence of sets $B'_{1} \subseteq B'_{2} \subseteq \cdots \subseteq B'_{i}$ satisfying properties~\eqref{CS:prop1} and~\eqref{CS:prop2}.
For $j \geq 1$ we define the finite set
\[
B_{N^{(i)}_j}^{(i+1)} = B_{N^{(i)}_j} \cap 
 \left[-\frac{i+1+g^{\ast}}{v}, \frac{i+1+g^{\ast}}{v}\right].
\]
Then $\left\{ B_{N^{(i)}_j}^{(i+1)} \right\}_{j=1}^{\infty}$ is an infinite sequence of subsets of the finite set 
$[-(i+1+g^{\ast})/v , (i+1+g^{\ast})/v ] \cap \Z$.  
By the pigeonhole principle, there is a strictly increasing sequence 
$\left\{N^{(i+1)}_j \right\}_{j=1}^{\infty}$ of positive integers and a set $B'_{i+1}$ such that $i+1 \leq N^{(i+1)}_1$ and 
\[
B'_i \subseteq B'_{i+1} = B^{(i+1)}_{N^{(i+1)}_j} \subseteq B_{N^{(i+1)}_j}
\]
for all $j\geq 1.$  Properties~\eqref{CS:prop1} and~\eqref{CS:prop2} are satisfied for $i+1.$  This completes the induction.  Moreover, the  set $A_h = \bigcup_{i=1}^{\infty} B'_{i}$ satisfies $R^{(\varphi)}_{\mca,A_h}(n) = f(n)$ for all $n \in \Z.$   This completes the proof.
\epf

\bt  
Let $h \geq 1$ and $\varphi(x_1,\ldots, x_h,y) = u_1x_1 + \cdots + u_hx_h + y.$   Let $\mca = (A_1,\ldots, A_{h})$
be an $h$-tuple of nonempty finite sets of integers and let $t \geq 1.$ 
Suppose that there is  a strictly increasing sequence $\{L_N\}_{N=1}^{\infty}$ of positive integers such that, for every $N \geq 1$, there exists a set $B_N$ of integers and a set $I_N$ consisting of $2L_N+1$ consecutive integers such that 
\[
R_{\mca,B_N}(n) = t \qquad\text{for $n \in I_N.$}
\]
Then there exists a set $B$ such that 
\[
R_{\mca,B}(n) = t \qquad\text{for all $n \in \Z.$}
\]
\et

\bpf
For every integer $N \geq 1$, there is an integer $c_N$ such that  $I_N = [c_N-L_N,c_N+L_N]\cap \Z.$  Replace the set $B_N$ with the set $B_N - c_N$ and apply Theorem~\ref{CS:theorem:compact}.
This completes the proof.
\epf

A related result appears in Nathanson~\cite{nath04c}.

\providecommand{\bysame}{\leavevmode\hbox to3em{\hrulefill}\thinspace}
\providecommand{\MR}{\relax\ifhmode\unskip\space\fi MR }
\providecommand{\MRhref}[2]{%
  \href{http://www.ams.org/mathscinet-getitem?mr=#1}{#2}
}
\providecommand{\href}[2]{#2}

\end{document}